\documentclass[11pt]{article}

\usepackage[T2A]{fontenc}
\usepackage[cp1251]{inputenc}
\usepackage{amsthm}
\usepackage{amsfonts}
\usepackage{eufrak}
\usepackage{amssymb}
\usepackage{amsmath}
\usepackage{cite}


\begin{document}
\newtheoremstyle{mytheorem}
  {\topsep}   
  {\topsep}   
  {\itshape}  
  {}       
  {\bfseries} 
  {  }         
  {5pt plus 1pt minus 1pt} 
  { }          
\newtheoremstyle{myremark}
  {\topsep}   
  {\topsep}   
  {\upshape}  
  {}       
  {\bfseries} 
  {  }         
  {5pt plus 1pt minus 1pt} 
  { }          
\theoremstyle{mytheorem}
\newtheorem{theorem}{Theorem}[section]
 \newtheorem*{theorema}{Theorem}
 \newtheorem*{A*}{Theorem A}
 \newtheorem*{B*}{Theorem B}
 \newtheorem*{heyde*}{The Heyde theorem}
 \newtheorem*{a*}{The Skitovich--Darmois theorem}
 \newtheorem*{b*}{The Skitovich--Darmois theorem for the group $\mathbb{R}\times \mathbb{Z}(2)$}
\newtheorem{proposition}[theorem]{Proposition}
\newtheorem{lemma}[theorem]{Lemma}
\newtheorem{corollary}[theorem]{Corollary}
\newtheorem{definition}[theorem]{Definition}
\theoremstyle{myremark}
\newtheorem{remark}[theorem]{Remark}

\noindent This article with minor changes was accepted for publication in 

\noindent Journal of Fourier Analysis and Applications

\bigskip

\noindent{\textbf{ The Heyde characterization theorem on compact totally }}

\noindent{\textbf{   disconnected and connected Abelian groups}}

\bigskip

\noindent{\textbf{Gennadiy Feldman}}

\bigskip

\noindent{\bf Abstract}

\bigskip

\noindent By the well-known  Heyde  theorem,  the Gaussian distribution
on the real line is characterized by the symmetry of the conditional
distribution of one linear form of  independent random variables
 given another. In the case of two independent random variables we give a complete
 description of compact totally disconnected Abelian groups $X$, where an analogue of  this theorem
 is valid.  We also  prove   that even a weak analogue of the Heyde theorem fails on compact connected
Abelian groups $X$. Coefficients of considered linear forms are topological automorphisms of $X$. The proofs are based on the study of solutions of a functional equation  on the character group of the group $X$ in the class of Fourier transforms of probability distributions.

\bigskip

\noindent{\bf Keywords}  Compact totally disconnected Abelian group $\cdot$ Compact connected Abelian group $\cdot$ Topological automorphism $\cdot$
  Haar  distribution
\bigskip

\noindent{\bf Mathematics Subject Classification} Primary  43A35 Secondary 60B15  $\cdot$ 62E10

\bigskip

\section{  Introduction}

Let $(\Omega, \mathfrak{A}, P)$ be a probability space,   $X$ be a topological Abelian group, $\mathfrak{B}$ be the     $\sigma$-algebra of Borel sets in $X$. Consider   random variables $\xi$ and $\eta$ defined on    $(\Omega, \mathfrak{A}, P)$ with values in   $(X, \mathfrak{B})$. The definition of the conditional distribution of  the random variable $\eta$ given $\xi$ see e.g. in  \cite[Appendix, \S 5]{B}. We note that the conditional distribution of  $\eta$ given $\xi$ is symmetric if and only if the random vectors $(\xi, \eta)$ and $(\xi, -\eta)$ are identically distributed.

According to the classical Skitovich--Darmois theorem the independence  of two linear forms of
  independent random variables is a characteristic property of the Gaussian distribution
on the real line. Another important characterization
result was proved by C.C. Heyde, where
instead of the independence the symmetry of the conditional distribution
of one linear form  given another is considered (\!\!\cite{He}, see
also \cite[\S\,13.4.1]{KaLiRa}). For two independent random variables
    the Heyde  theorem can be formulated as follows.
\begin{A*} \label{A*}
Let $\xi_1$  and $\xi_2$ be independent real random variables with distributions
       $\mu_1$ and $\mu_2$. Let $\alpha_j, \beta_j$ be non zero real numbers such that
$\alpha_1\beta_1^{-1}+\alpha_2\beta_2^{-1}\not=  0$.
If the conditional distribution of  the linear form  $L_2 = \beta_1\xi_1 +
\beta_2\xi_2$  given $L_1 = \alpha_1\xi_1 + \alpha_2\xi_2$ is symmetric, then
  $\mu_j$   are Gaussian distributions.
\end{A*}
Some generalisations of the Heyde  theorem, where independent
random variables
take values in a locally compact Abelian group $X$, and coefficients
of linear forms are topological automorphisms of $X$, were studied
in   \cite{Fe2, Fe4, Fe3, Fe20bb, Fe6, F, Fe7, {Fe8a}, Fe2020, POTA, My2, My1, MyF}, see also \cite[Chapter VI]{Fe5a}. We  continue  this research. The main result of the article is a complete description of compact totally disconnected
 Abelian groups for which a natural analogue of the Heyde  theorem for two independent random variables holds true.  It is well known that Gaussian distributions on a locally compact totally disconnected
 Abelian group $X$ are
degenerated (\!\!\cite[Chapter IV, \S 6, Remark 2]{Pa}). The role of  Gaussian distributions on such groups  is played by
shifts of the Haar distributions
of  compact subgroups  of $X$. Thus, we will be interested in the description of compact totally disconnected
 Abelian groups $X$ for which   shifts of the Haar distributions
of  compact subgroups  of $X$ are characterized by the symmetry of the conditional distribution
of one linear form of two independent random variables with values in $X$
 given another. In  \S 2
we give   a complete
 description of such groups. In \S 3 we prove a generalization of this result. In \S 4 we prove that even a weak analogue of the Heyde  theorem fails on compact connected
Abelian groups. The proofs are based on the study of solutions of a functional equation  on the character group of the group $X$ in the class of Fourier transforms of probability distributions.

In the article we  use standard results of abstract harmonic analysis, in particular  the duality theory for locally compact Abelian groups (see \cite{Hewitt-Ross}). Let $X$ be a second countable locally compact Abelian group. We
consider only such groups
without mentioning it specifically. Denote by $Y$ the character group
of the group $X$, and by  $(x,y)$ the value of a character $y \in Y$ at an
element $x \in X$. If $K$ is a closed subgroup of $X$, denote by
 $A(Y, K) = \{y \in Y: (x, y) = 1$ \mbox{ for all } $x \in K \}$
its annihilator.   Let $X_1$ and
$X_2$ be locally compact Abelian groups with character groups
$Y_1$ and $Y_2$ respectively. Let
 $\alpha:X_1\rightarrow X_2$
be a continuous homomorphism. The adjoint homomorphism $\tilde\alpha: Y_2\rightarrow Y_1$
is defined by the formula $(\alpha x_1,
y_2)=(x_1 , \tilde\alpha y_2)$ for all $x_1\in X_1, \ y_2\in
Y_2$. Denote by ${\rm Aut}(X)$ the group of topological automorphisms of
$X$ and  by $I$ the identity
automorphism of a group.

 Let
 $\{G_\iota: \iota\in {\cal I}\}$ be a nonempty set of compact Abelian
 groups. Denote by  $\mathop{\mbox{\rm\bf P}}\limits_{\iota \in
{\cal I}}G_\iota$ the  direct product of the groups
 $G_\iota$  equipped with the product topology.
 If $\{H_\iota: \iota\in {\cal I}\}$ is a nonempty set of discrete Abelian
 groups, then denote by  $\mathop{\mbox{\rm\bf P}^*}\limits_{\iota \in
{\cal I}}H_\iota$ the weak direct product of the groups
 $H_\iota$  equipped with   the discrete topology. We will always consider a finite Abelian group in the discrete topology.
Denote by   $\mathbb{Z}(n)=\{0, 1, \dots, n-1\}$ the group of  residue  classes modulo $n$.
Let $p$ be a prime number.
Denote by $\Delta_p$  the group of $p$-adic integers. A $p$-group is defined to be a group whose elements have order equal to a non-negative
power of $p$. In particular, for power $0$, it is possible that $X_p=\{0\}$.
Denote by $\cal P$
the set of prime numbers.

Denote by ${\rm M}^1(X)$ the
convolution semigroup of probability distributions on the group $X$.
 Let  $\mu\in{\rm M}^1(X)$. Denote by  $$\hat\mu(y) =
\int_{X}(x, y)d \mu(x), \ \ y\in Y,$$  the characteristic function (Fourier transform) of
the distribution  $\mu$, and by $\sigma(\mu)$ the support of $\mu$. Define the distribution $\bar \mu \in {\rm M}^1(X)$ by the formula
 $\bar \mu(B) = \mu(-B)$ for any Borel subset $B$ of $X$.
Then $\hat{\bar{\mu}}(y)=\overline{\hat\mu(y)}$.

Denote by $m_K$ the Haar distribution of a compact subgroup
 $K$ of the group $X$, and by  $I(X)$ the set of  shifts of Haar distributions $m_K$  of compact
subgroups
  $K$ of the group $X.$
We note that the characteristic function of a distribution
$m_K$ is of the form
\begin{equation}\label{11a}
\hat m_K(y)=
\begin{cases}
1, & \text{\ if\ }\   y\in A(Y, K),
\\  0, & \text{\ if\ }\ y\not\in
A(Y, K).
\end{cases}
\end{equation}
A distribution  $\gamma\in {\rm M}^1(X)$  is called Gaussian
(\!\!\cite[Chapter IV, \S 6]{Pa}, see also \cite{Fe1979})
if its characteristic function is represented in the form
$$
\hat\gamma(y)= (x,y)\exp\{-\varphi(y)\},
$$
where $x \in X$, and $\varphi(y)$ is a continuous nonnegative function
on the group $Y$
 satisfying the equation
 $$
\varphi(u + v) + \varphi(u
- v) = 2[\varphi(u) + \varphi(v)], \ \ u,  v \in
Y.
$$
Denote by $E_x$  the degenerate distribution
 concentrated at an element $x\in X$. Observe that degenerate distributions are Gaussian.
 Denote by $\Gamma(X)$ the set of Gaussian distributions on
    $X$.

\section{  The Heyde  theorem on compact totally disconnected
Abelian groups}

Let $\xi_j$, $j=1, 2$,   be   independent random variables with values in a locally compact Abelian group $X$ and distributions $\mu_j$, and let $\alpha_j, \beta_j\in {\rm Aut}(X)$. Assume that the conditional distribution of  the linear form  $L_2 = \beta_1\xi_1 +
\beta_2\xi_2$  given $L_1 = \alpha_1\xi_1 + \alpha_2\xi_2$ is symmetric. Obviously, if we are interested in the description of $\mu_j$, then we can suppose without loss of generality that $L_1 = \xi_1 + \xi_2$ and $L_2 = \xi_1 +
\alpha\xi_2$, where $\alpha\in {\rm Aut}(X)$.

Denote by $G$ the subgroup of $X$
generated by all elements of $X$ of order 2. Then $\alpha(G)=G$ and $g=-g$ for all $g\in G$.  Assume that $\xi_1$ and  $\xi_2$ take values in
$G$. Taking into account that the conditional distribution of  the linear form  $L_2 = \xi_1 +
\alpha\xi_2$  given $L_1 = \xi_1 + \xi_2$ is symmetric if and only if the random vectors $(L_1, L_2)$ and $(L_1, -L_2)$ are identically distributed, we see that   for any $\mu_j$ and any
topological automorphism $\alpha$ the conditional distribution of    $L_2$  given $L_1$ is symmetric.

Put  $K={\rm Ker}(I+\alpha)$ and suppose that $K\ne\{0\}$. Then $\alpha(k)=-k$ for all $k\in K$.
 Assume that
  $\xi_1$ and
$\xi_2$ are arbitrary  independent identically distributed random variables with values in
  $K$. Consider the linear form $L_1 = \xi_1 + \xi_2$ and $L_2 = \xi_1 +
\alpha\xi_2$.  Then $L_2 = \xi_1 -\xi_2$, and the characteristic functions of the random vectors $(L_1, L_2)$ and $(L_1, -L_2)$ are equal. Hence  the random vectors $(L_1, L_2)$ and $(L_1, -L_2)$ are identically distributed. So,
the conditional distribution of the linear form
    $L_2 = \xi_1 +
\alpha\xi_2$  given $L_1 = \xi_1 + \xi_2$  is symmetric.

Let $X$ be a compact totally disconnected
 Abelian group. Our goal is to describe all groups $X$ for which   shifts of the Haar distributions of  compact subgroups  of $X$ are characterized by the symmetry of the conditional distribution
of the linear form
    $L_2 = \xi_1 +
\alpha\xi_2$  given $L_1 = \xi_1 + \xi_2$. From what has been said above it follows that we should suppose that $X$ contains no elements of order 2, and
  condition
\begin{equation}\label{1}
{\rm Ker}(I+\alpha)=\{0\}
\end{equation}
is satisfied.

The main result of this section is the following theorem.
\begin{theorem}\label{th1}Let $X$ be a compact totally disconnected
Abelian group containing  no elements of order $2$, and let  $\alpha$ be a topological automorphism of
the group $X$ satisfying   condition $(\ref{1})$. Let
  $\xi_1$ and  $\xi_2$ be independent random variables with values in
       $X$  and distributions $\mu_1$ and $\mu_2$.
The symmetry of the   conditional distribution of the linear form
$L_2 = \xi_1 + \alpha\xi_2$ given  $L_1 = \xi_1 +
\xi_2$  implies that $\mu_j\in I(X)$, $j=1, 2$, if and only if the group $X$  is topologically isomorphic to a group of the form
\begin{equation}\label{n9}
 \mathop{\mbox{\rm\bf P}}\limits_{p\in {\cal P}, \ p>2}{X}_{p},
\end{equation}
where $X_p$ is a finite $p$-group.
\end{theorem}

Since shifts of the Haar distributions
of  compact subgroups  of the group $X$  play the role of Gaussian distributions on
$X$, Theorem \ref{th1} gives a complete
description of compact totally disconnected Abelian groups where an analogue of  the Heyde  theorem for two independent random variables holds. We note that an analogue of  the Heyde  theorem for two independent random variables is valid on an arbitrary discrete Abelian group (see \cite[Theorem 1]{Fe8a}).

To prove Theorem \ref{th1} we need some lemmas.
\begin{lemma} [\!\!{\protect\cite[Lemma 16.1]{Fe5a}}] \label{le1}   Let $X$ be a locally compact Abelian group, and let $\alpha$ be a topological automorphism of the group $X$.
 Let $\xi_1$ and $\xi_2$ be independent random variables with values in
       $X$  and distributions
 $\mu_1$ and $\mu_2$.
 The conditional distribution of the linear form
 $L_2 = \xi_1 + \alpha\xi_2$
 given $L_1 = \xi_1 + \xi_2$ is symmetric if and only
 if the characteristic functions
 $\hat\mu_j(y)$ satisfy the equation
\begin{equation}\label{2a}
\hat\mu_1(u+v )\hat\mu_2(u+\tilde\alpha v )=
\hat\mu_1(u-v )\hat\mu_2(u-\tilde\alpha v), \ \ u, v \in Y.
\end{equation}
\end{lemma}
\begin{lemma}    [\!\!{\protect\cite {Fe2}, see also \cite[Corollary 17.2 and Remark 17.5]{Fe5a}}]
\label{le2} Let  $X$ be a  finite Abelian group
  containing no elements of order $2$. Let  $\alpha$ be an automorphism of
 the group $X$ satisfying   condition $(\ref{1})$.
Let
  $\xi_1$ and  $\xi_2$ be independent random variables with values in
       $X$  and distributions $\mu_1$ and $\mu_2$.
If the conditional distribution of the linear form
$L_2 = \xi_1 + \alpha\xi_2$ given  $L_1 = \xi_1 +
\xi_2$  is symmetric, then
$\mu_j=m_K*E_{x_j}$, where $K$ is a
subgroup of  $X$,
$x_j\in X$, $j=1, 2$. Moreover, $\alpha (K)=K$.
\end{lemma}
\begin{lemma}[\!\!{\protect\cite[Theorem 2] {F}}]\label{le4} Let $X=\Delta_p$ be the group of $p$-adic integers.
Then there exist a topological automorphism $\alpha$ of
the group $X$ satisfying   condition $(\ref{1})$, and
   independent random variables $\xi_1$ and  $\xi_2$  with values in
       $X$  and distributions $\mu_1$ and $\mu_2$ such that
 the conditional distribution of the linear form
$L_2 = \xi_1 + \alpha\xi_2$ given  $L_1 = \xi_1 +
\xi_2$  is symmetric, while  $\mu_j\notin
I(X)$, $j=1, 2$.
\end{lemma}
The following lemma plays a key role in the proof of Theorem \ref{th1}.
\begin{lemma}  \label{le5} Let $p$ be a prime number, and let $p>2$. Consider the group
\begin{equation}\label{n1}
X=
\mathop{\mbox{\rm\bf P}}\limits_{n=1}^\infty \mathbb{Z}(p^{k_n}),
\ \   k_n\le k_{n+1}, \ \ n=1, 2, \dots
\end{equation}
Then there exist a topological automorphism $\alpha$ of
the group $X$ satisfying   condition $(\ref{1})$, and
   independent identically distributed
   random variables $\xi_1$ and  $\xi_2$  with values in
       $X$  and distribution $\mu$   such that
 the conditional distribution of the linear form
$L_2 = \xi_1 + \alpha\xi_2$ given  $L_1 = \xi_1 +
\xi_2$  is symmetric, while $\mu\notin
I(X)$.
\end{lemma}
\textbf{\textit{Proof }}   The group $Y$
is topologically isomorphic to the weak direct product of the groups
 $\mathbb{Z}(p^{k_n})$.  In order not to complicate the notation we assume that
\begin{equation}\label{n2}
Y=
{\mathop{\mbox{\rm\bf P}^*}\limits_{n=1}^\infty}
 \mathbb{Z}(p^{k_n}).
\end{equation}
We denote by $t=\{t_n\}_{n=1}^\infty$, $t_n\in\mathbb{Z}(p^{k_n})$, elements
of the group
 $X$ and by  $s=\{s_n\}_{n=1}^\infty$, $s_n\in\mathbb{Z}(p^{k_n})$, elements
 of the group $Y$. Let $s=\{s_n\}_{n=1}^\infty\in Y$. We say that $s\in \mathbb{Z}(p^{k_1})$ if $s_2=s_3=\dots=0$.

 Let $i\le j$. By $\alpha_{i, j}$ denote
the homomorphism $\alpha_{i, j}:\mathbb{Z}(p^{k_i})\rightarrow\mathbb{Z}(p^{k_j})$
of the form
$$
\alpha_{i, j}t_i=p^{k_j-k_i}t_i, \ \ t_i\in \mathbb{Z}(p^{k_i}).
$$
It is easy to verify that the adjoint homomorphism $\tilde\alpha_{i, j}:\mathbb{Z}(p^{k_j})\rightarrow\mathbb{Z}(p^{k_i})$ is of the form
$$
\tilde\alpha_{i, j}s_j=s_j ({\rm mod} \ p^{k_i}), \ \ s_j\in \mathbb{Z}(p^{k_j}).
$$
Define a continuous homomorphism $\alpha:X\rightarrow X$ by the formula $\alpha\{t_n\}_{n=1}^\infty=\{h_n\}_{n=1}^\infty$, where
\begin{equation}\label{n3}
h_n=
\begin{cases}
(p^{k_1}-t_1) ({\rm mod} \ p^{k_1}), & \text{\ if\ }\   n=1,
\\  (\alpha_{n-1, n}t_{n-1}+p^{k_n}-t_n)({\rm mod} \ p^{k_n}), & \text{\ if\ }\ n\ge 2.
\end{cases}
\end{equation}
\bigskip
It is easy to see that $\alpha$ is a topological automorphism
of the group
$X$ satisfying   condition (\ref{1}). We find from (\ref{n3}) that
the adjoint automorphism $\tilde\alpha:Y\rightarrow Y$ is of the form
$\tilde\alpha\{s_n\}_{n=1}^\infty=\{u_n\}_{n=1}^\infty$, where
\begin{equation}\label{n4}
u_n=(\tilde\alpha_{n, n+1}s_{n+1}+p^{k_n}-s_n)({\rm mod} \ p^{k_n}), \ \ n\ge 1.
\end{equation}

It follows from (\ref{n1})  and (\ref{n2}) that  $X=\mathbb{Z}(p^{k_1})\times G$, where
$G=\mathop{\mbox{\rm\bf P}}\limits_{n=2}^\infty \mathbb{Z}(p^{k_n})$ and
  $Y=\mathbb{Z}(p^{k_1})\times H$, where
$H=\mathop{\mbox{\rm\bf P}^*}\limits_{n=2}^\infty \mathbb{Z}(p^{k_n})$.
Obviously,  $A(Y, G)=\mathbb{Z}(p^{k_1})$.

Take $0<a<1$. Consider on the group $X$
the distribution $\mu=a m_{X}+(1-a)m_{G}.$
  It follows from (\ref{11a}) that the characteristic function
 $\hat\mu(y)$ is of the form
\begin{equation}\label{n5}
\hat\mu(y)=
\begin{cases}
1, & \text{\ if\ }\   y=0,
\\  1-a, & \text{\ if\ }\ y \in
{\mathbb{Z}}(p^{k_1})\backslash\{0\},
\\ 0, & \text{\ if\ }\ y \notin
{\mathbb{Z}}(p^{k_1}).
\end{cases}
\end{equation}
Let $\xi_1$ and  $\xi_2$ be independent identically distributed
   random variables   with values in
       $X$  and distribution $\mu$.  Let us verify that
 the conditional distribution of the linear form
$L_2 = \xi_1 + \alpha\xi_2$ given  $L_1 = \xi_1 +
\xi_2$  is symmetric.
By Lemma \ref{le1}, it suffices to verify that the characteristic function $\hat\mu(y)$
satisfies the equation
  \begin{equation}\label{08a2}
\hat\mu(u+v )\hat\mu(u+\tilde\alpha v )=
\hat\mu(u-v )\hat\mu(u-\tilde\alpha v), \ \ u, v \in Y.
\end{equation}

Let $u, v \in Y$. Consider 3 cases.

1. $u,  v \in {\mathbb{Z}}(p^{k_1})$. It follows from (\ref{n4}) that
 $\tilde\alpha y = -y$ when  $y \in {\mathbb{Z}}(p^{k_1})$. Hence, the restriction
of equation (\ref{08a2}) to the subgroup
  ${\mathbb{Z}}(p^{k_1})$ takes the form
$$
 \hat\mu(u+v)  \hat\mu(u-v)=  \hat\mu(u-v)   \hat\mu(u+v), \ \ u, v
\in {\mathbb{Z}}(p^{k_1}).
$$
It is obvious that (\ref{08a2}) is fulfilled.

2. Either $u \in {\mathbb{Z}}(p^{k_1}),$ $v \notin {\mathbb{Z}}(p^{k_1})$  or
$u \notin {\mathbb{Z}}(p^{k_1}),$ $v \in {\mathbb{Z}}(p^{k_1}).$ Then  $u
\pm v \notin {\mathbb{Z}}(p^{k_1})$. It follows from (\ref{n5}) that $\hat\mu(u \pm v)=0$.
Then both sides of equation (\ref{08a2}) are equal to zero.

 3.  $u, v \notin {\mathbb{Z}}(p^{k_1}).$ Assume that the left-hand side in
  (\ref{08a2}) is not equal to zero. Then (\ref{n5}) implies that
  $u+v, u+\tilde\alpha v \in
{\mathbb{Z}}(p^{k_1})$. It follows from this that
\begin{equation}\label{n7}
  (I-\tilde\alpha )v
\in {\mathbb{Z}}(p^{k_1}).
\end{equation}
Let $v=\{v_n\}_{n=1}^\infty$, $v_n\in\mathbb{Z}(p^{k_n})$. We find from (\ref{n4}) and (\ref{n7}) that
\begin{equation}\label{n8}
(2v_n+p^{k_n}-\tilde\alpha_{n, n+1}v_{n+1})({\rm mod} \ p^{k_n})=0, \ n\ge 2.
\end{equation}
It follows from  $p>2$ that if $v_n\ne 0$,  then $2v_n({\rm mod} \ p^{k_n})\ne 0$, $n\ge 1$. Assume that
  $v_2\ne 0$. Then we  find from (\ref{n8}) that $v_n\ne 0$ for each $n\ge 2$.
  But this contradicts the fact that  $v_n=0$ for all but a finite set of indices $n$.   Hence, $v_2=0$. Reasoning by induction we   prove  that
  $v_3=0$, then $v_4= 0$, and  so on.
Thus, we got that $v=\{v_n\}_{n=1}^\infty\in {\mathbb{Z}}(p^{k_1})$, contrary to assumption.
The obtained contradiction shows that the left-hand side of
  (\ref{08a2}) is  equal to zero. Similarly, we prove that when $u, v \notin {\mathbb{Z}}(p^{k_1})$ the right-hand side of
  (\ref{08a2}) is also equal to zero. So,    both sides of equation (\ref{08a2}) are equal to zero.
   $\blacksquare$

   The following lemma follows directly from Lemma 13.24 in \cite{Fe5a}.
\begin{lemma} \label{le6} Let  $X$ be a compact
totally disconnected Abelian group containing no elements of
order $2$.
Then either    $X$ is
topologically isomorphic to a group of the form $(\ref{n9})$
or for some prime number   $p$ there exists a compact subgroup $K$ of $X$ such that $K$ is
a topological direct factor  of the group  $X$, and $K$ is
topologically isomorphic to either the group of $p$-adic integers $\Delta_p$  or a group of the form $(\ref{n1})$.
\end{lemma}
\textbf{\textit{Proof of Theorem  \ref{th1}}}   Necessity. Assume that a compact totally disconnected
Abelian group $X$ contains no elements of order $2$, and
$X$ is not topologically isomorphic to a group of the form $(\ref{n9})$.
 By Lemma \ref{le6},  for some prime number   $p$, $p>2$,  there exists a compact subgroup $K$ of the group $X$ such that $K$ is
a topological direct factor  of    $X$, and $K$ is
topologically isomorphic to either a group of $p$-adic integers $\Delta_p$  or a group of the form $(\ref{n1})$. We have $X=K\times G$. Denote elements of the group $X$ by $x=(k, g)$, where $k\in K$, $g\in G$. Applying either Lemma \ref{le4}  or  Lemma \ref{le5} we find a topological automorphism $\alpha_K$ of
the group $K$ satisfying   condition $(\ref{1})$, and independent random variables $\xi_1$ and
$\xi_2$  with values in $K$ and distributions $\mu_1$ and $\mu_2$
such that the conditional distribution of the linear form $M_2=\xi_1
+ \alpha_K \xi_2$ given $M_1=\xi_1+\xi_2$ is symmetric, whereas $\mu_j\notin I(K)$, $j=1, 2$. Put
$\alpha(k, g)=(\alpha_K k, g)$. Since the group $X$ contains no
elements of order 2,  $\alpha$ is a topological automorphism   of
$X$ satisfying   condition $(\ref{1})$. Obviously, if we consider $\xi_j$ as independent random variables    with values in the group $X$, then the conditional distribution of the linear form $L_2=\xi_1
+ \alpha  \xi_2$ given $L_1=\xi_1+\xi_2$ is symmetric, whereas $\mu_j\notin I(X)$, $j=1, 2$. Necessity is proved.

  Sufficiency. Assume that a compact totally disconnected
Abelian group $X$  is topologically isomorphic to a group of the form (\ref{n9}), and let    $\alpha$ be a topological automorphism of
the group $X$ satisfying   condition $(\ref{1})$. Let
  $\xi_1$ and  $\xi_2$ be independent random variables with values in
       $X$  and distributions $\mu_1$ and $\mu_2$.
By Lemma \ref{le1}, the symmetry of the conditional distribution
of the linear form $L_2$ given $L_1$ implies that the characteristic functions $\hat\mu_j(y)$ satisfy equation
(\ref{2a}). Put
$\nu_j = \mu_j
* \bar \mu_j$. Then
 $\hat \nu_j(y) = |\hat \mu_j(y)|^2 \ge 0$  for all $y \in Y$,
 and  the characteristic functions  $\hat \nu_j(y)$
also satisfy equation (\ref{2a}).

It follows from  (\ref{n9}) that
\begin{equation}\label{m1}
 Y=\mathop{\mbox{\rm\bf P}^*}\limits_{p\in {\cal P}, \ p>2}{Y}_{p},
\end{equation}
where the group $Y_p$ is   isomorphic to the character group
of the group $X_p$.
Put
\begin{equation} \label{25.03n1}
A_m=\mathop{\mbox{\rm\bf P}}\limits_{p\in {\cal P}, \ 2<p\le m}{X}_{p}, \ \ B_m=\mathop{\mbox{\rm\bf P}}\limits_{p\in {\cal P}, \ 2<p\le m}{Y}_{p}, \ \ m=3, 4, \dots
\end{equation}
Then $A_m$ and $B_m$ are finite Abelian groups, $B_m$ is isomorphic to the character group of the group $A_m$, $B_m\subset  B_{m+1}$  and $Y=\bigcup\limits_{m=3}^\infty B_m$.
Denote by $\alpha_m$ the restriction of $\alpha$ to the subgroup $A_m$. It is obvious that $\alpha_m\in {\rm Aut}(A_m)$ and $\alpha_m$ satisfies the condition  $(\ref{1})$. Note that the restriction of $\tilde\alpha$ to   subgroup $B_m$ coincides with $\tilde\alpha_m$.
Consider the restriction of equation (\ref{2a}) for the characteristic functions $\hat\nu_j(y)$ to  the subgroup $B_m$. We have
\begin{equation}\label{25.03.1}
\hat\nu_1(u+v )\hat\nu_2(u+\tilde\alpha_m v )=
\hat\nu_1(u-v )\hat\nu_2(u-\tilde\alpha_m v), \ \ u, v \in B_m.
\end{equation}
Taking into account Lemma \ref{le1} and applying
Lemma \ref{le2} to the finite group $A_m$, we obtain from   (\ref{25.03.1}) that $\hat\nu_1(y)=\hat\nu_2(y)$ for all
$y\in B_m$, and the characteristic functions $\hat\nu_j(y)$
take only values 0 and 1 on $B_m$.  This implies that $\hat\nu_1(y)=\hat\nu_2(y)$ for all
$y\in Y$, and the characteristic  functions $\hat\nu_j(y)$
take only values 0 and 1 on $Y$. Put $E=\{y\in Y: \hat\nu_j(y)=1\}$. Then $E$  is a subgroup of $Y$. Put $K=A(X, E)$. Since $E=A(Y, K)$, it follows from (\ref{11a}) that $\nu_1=\nu_2=m_K$. This easily implies that
$\mu_j=m_K*E_{x_j}$, where $x_j\in X$, $j=1, 2$.
Sufficiency is proved.

In fact, we have proved somewhat more, namely $\mu_j$ are shifts of the same Haar distribution   of a compact subgroup $K$ of the group $X$. Furthermore, we can assert  that $\alpha(K)=K$. Indeed, put $E_m=E\cap B_m$.
It follows from Lemma \ref{le2} that $\tilde\alpha_m(E_m)=E_m$. Hence, $\tilde\alpha(E)=E$, so that $\alpha(K)=K$.
  $\blacksquare$

\begin{remark}\label{r1} Let us compare Theorem \ref{th1} with the Skitovich--Darmois theorem for compact totally disconnected
Abelian groups, which has the following statement. It can be deduced from   Theorem 1 in
\cite{FeGra1}, see also \cite[Theorem 13.25]{Fe5a}.
\begin{theorema} Let $X$ be a compact totally disconnected
Abelian group containing no elements of order $2$,  and let  $\alpha$ be a topological automorphism  of the group  $X$. Let
  $\xi_1$ and  $\xi_2$ be independent random variables with values in
       $X$  and distributions $\mu_1$ and $\mu_2$.
The independence of   the     linear forms $L_1 = \xi_1 +
\xi_2$ and
$L_2 = \xi_1 + \alpha\xi_2$   implies that $\mu_j\in I(X)$, $j=1, 2$, if and only if the group $X$  is topologically isomorphic to a group of the form
$$
 \mathop{\mbox{\rm\bf P}}\limits_{p\in {\cal P}}(\Delta_p^{n_p}\times {X}_{p}),
$$
where $n_p$ is a nonnegative integer, and  $X_p$ is a finite $p$-group, $X_2=\{0\}$.
\end{theorema}
We see that the class of compact totally disconnected Abelian groups, for which the Skitovich--Darmois theorem is true, is wider than the class of compact totally disconnected Abelian groups for which the Heyde  theorem is valid.

At the same time, we note that both the Skitovich--Darmois theorem and the Heyde theorem are valid for discrete Abelian groups   containing no elements of order 2.
\end{remark}

 We supplement Theorem \ref{th1} with the following statements.
\begin{proposition} \label{pr1} Let $X$ be a compact totally disconnected
Abelian group of the form $(\ref{n9})$, let $K$ be a compact subgroup of $X$, and let
  $\alpha$ be a topological automorphism of the group
 $X$ satisfying   condition $(\ref{1})$.
Let
  $\xi_1$ and  $\xi_2$ be independent identically distributed
  random variables with values in
       $X$  and distribution $m_K$. Then the following assertions  are equivalent:

$(i)$ the conditional distribution of the linear form
$L_2 = \xi_1 + \alpha\xi_2$ given  $L_1 = \xi_1 +
\xi_2$  is symmetric;

$(ii)$ $(I-\alpha )(K)=K$.
\end{proposition}
\textbf{\textit{Proof }}  $(i) \Rightarrow (ii)$. Note that every closed subgroup of a group of the form $(\ref{n9})$ is topologically isomorphic to a group of the form $(\ref{n9})$. By Theorem \ref{th1}, $\alpha(K)=K$, i.e. the restriction of $\alpha$ to the subgroup $K$ is a topological automorphism of $K$.  Thus, we can prove the statement assuming that $K=X$.   We note that $(I-\alpha )(X)=X$ if and only if ${\rm Ker}(I-\tilde\alpha)=\{0\}$. Denote by $\tilde\alpha_m$ the same topological automorphism as in the proof of  sufficiency in Theorem \ref{th1}. We have ${\rm Ker}(I-\tilde\alpha)=\{0\}$ if and only if ${\rm Ker}(I-\tilde\alpha_m)=\{0\}$ for   $m=3, 4, \dots$. Taking into account Lemma \ref{le1}, the statement  ${\rm Ker}(I-\tilde\alpha_m)=\{0\}$ follows from the fact that the proposition is valid in the case when $X$ is
  a  finite Abelian group
  containing no elements of order $2$ (\!\!\cite[Proposition 1]{Fe8a}).

 $(ii) \Rightarrow (i)$. This statement holds   for  an arbitrary locally compact Abelian group  $X$    and a compact subgroup  $K$    of $X$ (\!\!\cite[Proposition 1]{Fe8a}). $\blacksquare$
\begin{proposition} \label{pr5}
Let $X$ be a compact totally disconnected
Abelian group of the form
$$
X=\mathop{\mbox{\rm\bf P}}\limits_{p\in {\cal P}}{X}_{p},
$$
where $X_p$ is a finite $p$-group. Put $$G=\mathop{\mbox{\rm\bf P}}\limits_{p\in {\cal P}, \ p>2}{X}_{p}.$$  Let  $\alpha$ be a topological automorphism of
the group $X$ satisfying   condition $(\ref{1})$. Let
  $\xi_1$ and  $\xi_2$ be independent random variables with values in
       $X$  and distributions $\mu_1$ and $\mu_2$.
If the conditional distribution of the linear form
$L_2 = \xi_1 + \alpha\xi_2$ given  $L_1 = \xi_1 +
\xi_2$  is symmetric, then
$\mu_j=\rho_j*m_K*E_{g_j}$, where   $\rho_j$ are some distributions on $X_2$, $K$ is a compact subgroup
of    $G$,
$g_j\in G$, $j=1, 2$.
\end{proposition}
To prove Proposition \ref{pr5} we need the following lemma  which is a special case of Theorem 3 proved in \cite{Fe8a}.
\begin{lemma}  \label{le11} Let $X$ be a finite
Abelian group, let $X_2$ be its $2$-component, and let $G$ be the subgroup of $X$ generated by all elements of odd order.   Let  $\alpha$ be an  automorphism of
the group $X$ satisfying   condition $(\ref{1})$. Let
  $\xi_1$ and  $\xi_2$ be independent random variables with values in
       $X$  and distributions $\mu_1$ and $\mu_2$.
If the conditional distribution of the linear form
$L_2 = \xi_1 + \alpha\xi_2$ given  $L_1 = \xi_1 +
\xi_2$  is symmetric, then
$\mu_j=\rho_j*m_K*E_{g_j}$, where $\sigma(\rho_j)\subset X_2$, $K$ is a subgroup
of    $G$,
$g_j\in G$, $j=1, 2$.
\end{lemma}
\textbf{\textit{Proof of Proposition \ref{pr5}}}  We argue  as in the proof of sufficiency    in Theorem \ref{th1}.  We have $X=X_2\times G$.
  Then $Y= L\times H$, where $L$ and $H$ are groups
topologically isomorphic to the character groups of the groups $X_2$ and $G$ respectively. Denote by
 $y=(l, h)$, where $l\in L$, $h\in H$, elements of the group $Y$. Taking into account Lemma \ref{le1}, Theorem \ref{th1} and the fact that the restriction of $\alpha$ to the subgroup  $G$ is a topological automorphism of  $G$ satisfying   condition (\ref{1}), it is easy to verify that Proposition \ref{pr5} will be proved if we prove that there exist functions $a_j(l)$, $l\in L$, and $b_j(h)$, $h\in H$, such that $$\hat\mu_j(l, h)=a_j(l)b_j(h), \ \ l \in L, \  \ h\in H, \ \ j=1, 2.$$
This implies that $a_j(l)$ are the characteristic functions of some distributions $\rho_j$ on $X_2$ and $b_j(h)$ are the characteristic functions of some shifts of a distribution $m_K$, where $K$ is a subgroup
of    $G$.

We have $Y=\mathop{\mbox{\rm\bf P}^*}\limits_{p\in {\cal P}}{Y}_{p}$,
where the group $Y_p$ is   isomorphic to the character group
of the group $X_p$.
Define the groups $A_m$ and $B_m$ by      (\ref{25.03n1}).
Then $A_m$ and $B_m$ are finite Abelian groups, the group $L\times B_m$ is isomorphic to the character group of the group $X_2\times A_m$, $B_m\subset B_{m+1}$  and $Y=L\times\bigcup\limits_{m=3}^\infty B_m$. Denote by $\alpha_m$ the restriction of $\alpha$ to the subgroup $X_2\times A_m$. It is obvious that $\alpha_m\in {\rm Aut}(X_2\times A_m)$ and $\alpha_m$ satisfies condition  $(\ref{1})$. Note that the restriction of $\tilde\alpha$ to the subgroup $L\times B_m$ coincides with $\tilde\alpha_m$.
By Lemma \ref{le1}, the symmetry of the conditional distribution of the linear form
$L_2$ given $L_1$ implies that the characteristic functions  $\hat\mu_j(l, h)$ satisfy  equation (\ref{2a}). Consider the restriction of equation (\ref{2a}) to the subgroup $L\times B_m$. We have
\begin{equation}\label{24.03.2}
\hat\mu_1(u+v )\hat\mu_2(u+\tilde\alpha_m v )=
\hat\mu_1(u-v )\hat\mu_2(u-\tilde\alpha_m v), \ \ u, v   \in L\times B_m.
\end{equation}
Taking into account Lemma \ref{le1} and applying Lemma   \ref{le11} to the finite group $X_2\times A_m$,    we obtain from     (\ref{24.03.2})   that there exist functions $a_{j}(l)$, $l\in L$, and $b_{j, m}(h)$, $h \in B_m$ such that $\hat\mu_j(l, h)=a_j(l)b_{j, m}(h)$ for all $(l, h)\in L\times B_m$. Moreover,   $b_{j, m+1}(h)=b_{j, m}(h)$ for all $h \in B_m$. Put $b_{j}(h)=b_{j, m}(h)$ for all $h \in B_m$. Then
  $\hat\mu_j(l, h)=a_j(l)b_j(h)$ for all $(l, h)\in Y$, $j=1, 2$. $\blacksquare$

\section{ A generalization of Theorem \ref{th1}}

The following theorem has been proved in  \cite[Theorem 2]{Fe8a}.
\begin{B*} \label{B*}  Let $X={\mathbb R}^n\times D$, where $n \ge 0,$ and $D$ is
a discrete Abelian group
  containing no elements of order $2$. Let  $\alpha$ be a topological automorphism of the group
 $X$ satisfying   condition $(\ref{1})$.
 Let
  $\xi_1$ and  $\xi_2$ be independent random variables with values in
       $X$  and distributions $\mu_1$ and $\mu_2$.
If the conditional distribution of the linear form
$L_2 = \xi_1 + \alpha\xi_2$ given  $L_1 = \xi_1 +
\xi_2$  is symmetric, then
$\mu_j=\gamma_j*m_K*E_{x_j}$, where $\gamma_j \in \Gamma({\mathbb R}^n)$,
 $K$ is a finite
subgroup of   $D$,
$x_j\in D$, $j=1, 2$. Moreover, $\alpha (K)=K$.
\end{B*}
In this section we shall prove a theorem which generalizes both Theorem \ref{th1} and
Theorem B and describe a wide class of locally compact Abelian groups where an analogue of   the Heyde  theorem holds.
\begin{theorem} \label{th2}  Let
\begin{equation}\label{09.09.15.1}
 X={\mathbb R}^n\times D\times G,
\end{equation}
where $n \ge 0,$  $D$ is
a discrete Abelian group
  containing no elements of order $2$, and
   $G$ is  a compact totally disconnected
Abelian group  of the form $(\ref{n9})$.
Let  $\alpha$ be a topological automorphism of the group
 $X$ satisfying   condition $(\ref{1})$.
Let
  $\xi_1$ and  $\xi_2$ be independent random variables with values in
       $X$  and distributions $\mu_1$ and $\mu_2$.
If the conditional distribution of the linear form
$L_2 = \xi_1 + \alpha\xi_2$ given  $L_1 = \xi_1 +
\xi_2$  is symmetric, then
$\mu_j=\gamma_j*m_K*E_{x_j}$, where $\gamma_j \in \Gamma({\mathbb R}^n)$,
 $K$   is a compact
subgroup of   $X$,
$x_j\in X$, $j=1, 2$. Moreover, $\alpha (K)=K$.
\end{theorem}
To prove Theorem \ref{th2} we need two lemmas. One of them is  the following well-known statement
 (see e.g. \cite[Proposition 2.13]{Fe5a}).
\begin{lemma}\label{le7} Let $X$ be a locally compact Abelian group,
and  let  $\mu\in{\rm
M}^1(X)$. Then the set  $E=\{y\in Y:\ \hat\mu(y)=1\}$ is a closed
subgroup of the group
   $Y$,   and  $\sigma(\mu)\subset A(X,E)$.
\end{lemma}
The following lemma generalizes  Theorem \ref{th1} and plays a key role in the proof of Theorem \ref{th2}.
\begin{lemma} \label{le8} Let $X={\mathbb R}^n\times G$, where $n \ge 0,$
and $G$ is a compact totally disconnected
Abelian group
   of the form $(\ref{n9})$. Let  $\alpha$ be a topological automorphism of the group
 $X$ satisfying   condition $(\ref{1})$.
Let
  $\xi_1$ and  $\xi_2$ be independent random variables with values in
       $X$  and distributions $\mu_1$ and $\mu_2$.
If the conditional distribution of the linear form
$L_2 = \xi_1 + \alpha\xi_2$ given  $L_1 = \xi_1 +
\xi_2$  is symmetric, then
$\mu_j=\gamma_j*m_K*E_{x_j}$, where $\gamma_j \in \Gamma({\mathbb R}^n)$,
  $K$ is a compact
subgroup of    $G$,
$x_j\in G$, $j=1, 2$. Moreover, $\alpha (K)=K$.
\end{lemma}
\textbf{\textit{Proof }}    Denote elements of the group $X$ by $x=(t, g)$, where
$t\in {\mathbb R}^n$, $g\in G$.
We have $Y={\mathbb R}^n\times H$, where $H$ is a group topologically isomorphic
to the character group of the group $G$. Denote elements of the group $Y$ by
 $y=(s, h)$, where $s\in {\mathbb R}^n$, $h\in H$.
Since ${\mathbb R}^n$ is the connected component of zero of
the group $X$, we have $\alpha({\mathbb R}^n)={\mathbb R}^n$. Since
 $G$ is the subgroup of $X$ consisting of all compact elements of the group
  $X$, we have $\alpha(G)=G$.
This implies that the restriction of
 $\alpha$ to each of the subgroups ${\mathbb R}^n$ and $G$ is a topological automorphism
of the corresponding  subgroup. Denote by $\alpha_{{\mathbb R}^n}$
and $\alpha_G$ these restrictions respectively.  This means that  $\alpha$ can be
written in the form  $\alpha(t, g)= (\alpha_{{\mathbb R}^n} t, \alpha_G g)$, $(t,
g) \in X$. Note that the restrictions of $\tilde\alpha$ to the subgroups ${\mathbb R}^n$ and   $H$ coincide  with $\tilde\alpha_{{\mathbb R}^n}$ and $\tilde\alpha_G$ respectively. Hence, $\tilde\alpha$ can be
written in the form   $\tilde\alpha(s, h)= (\tilde\alpha_{{\mathbb R}^n} s, \tilde\alpha_G h)$, $(s,
h) \in Y$.

By Lemma \ref{le1}, the symmetry of the conditional distribution of the linear form
$L_2$ given $L_1$ implies that the characteristic functions  $\hat\mu_j(s, h)$ satisfy
 equation
(\ref{2a}) which takes the form
\begin{equation}\label{04.09.1}
\hat\mu_1(s_1+s_2, h_1+h_2)\hat\mu_2(s_1+\tilde\alpha_{{\mathbb R}^n}s_2, h_1+\tilde\alpha_G h_2)$$$$=
\hat\mu_1(s_1-s_2, h_1-h_2)\hat\mu_2(s_1-\tilde\alpha_{{\mathbb R}^n}s_2, h_1-\tilde\alpha_G h_2), \ \  s_j \in {\mathbb R}^n, \ \ h_j \in H.
\end{equation}
It is obvious that the topological automorphism $\alpha_{{\mathbb R}^n}$
 satisfies   condition (\ref{1}). Substituting  $h_1=h_2=0$ in (\ref{04.09.1}),
 taking into account Lemma \ref{le1} and
 Theorem B for the group $X={\mathbb R}^n$, we obtain from the resulting equation
 that
\begin{equation}\label{04.09.2}
\hat\mu_j(s, 0)=\exp\{i\langle t_j, s\rangle-\langle A_js,
s\rangle\},\ \ s\in \mathbb{R}^n, \ \ j=1, 2,
\end{equation}
where $t_j\in \mathbb{R}^n$, $\langle ., .\rangle$ is the scalar product,
and $A_j$ are symmetric positive
semi-definite $n\times n$-matrices.
Similarly, substituting  $s_1=s_2=0$ in (\ref{04.09.1}), taking into account Lemma \ref{le1} and   applying
 Theorem \ref{th1}  to the group $X=G$, we obtain from the resulting equation
 that
\begin{equation}\label{04.09.3}
\hat\mu_j(0, h)=\hat m_K(h)(g_j, h), \ \ h\in H, \ \ j=1, 2,
\end{equation}
where $K$ is a compact subgroup in $G$, and $g_j\in G$. Moreover, $\alpha_G(K)=K$.

Substituting (\ref{04.09.3})  into (\ref{04.09.1})  and taking into account
 (\ref{11a}), we find that
$2(g_1+\alpha_Gg_2)\in K$. Since $G$ is a group
   of the form $(\ref{n9})$  and $K$ is a compact subgroup of $G$, multiplication by
   2 is a topological automorphism of $K$. Hence,
   \begin{equation}\label{03.04.1}
g_1+\alpha_Gg_2\in K.
\end{equation}
 Put  $x_1=-\alpha_Gg_2$, $x_2=g_2$. Then
 $x_1+\alpha_Gx_2=0$. It follows from (\ref{11a}),   (\ref{04.09.3}), and (\ref{03.04.1}) that
\begin{equation}\label{04.09.4}
\hat\mu_j(0, h)=\hat m_K(h)(x_j, h), \ \ h\in H, \ \ j=1, 2.
\end{equation}
Consider new independent random variables  $\eta_j=\xi_j-x_j$, and denote by $\lambda_j$ the distribution of the random variable
$\eta_j$, $j=1, 2$. Since
\begin{equation}\label{04.09.5}
\lambda_j=\mu_j*E_{-x_j}, \ \ j=1, 2,
\end{equation}
and $x_1+\alpha x_2=0$,   the characteristic functions
$\hat\lambda_j(y)$ also satisfy equation (\ref{2a}). By Lemma \ref{le1},
the conditional distribution of the linear form
$N_2 = \eta_1 + \alpha\eta_2$ given
  $N_1 = \eta_1 + \eta_2$   is symmetric. It follows from (\ref{11a}), (\ref{04.09.4}), and
(\ref{04.09.5}) that $\hat\lambda_j(0, h)=1$ for all $h\in A(H, K)$. By Lemma \ref{le7}, this implies that
 $\sigma(\lambda_j)\subset A(X, A(H, K))={\mathbb R}^n\times K$.
Since $\alpha({\mathbb R}^n\times K)={\mathbb R}^n\times K$, we can assume that the independent random variables
 $\eta_j$ take values in the group $X={\mathbb R}^n\times K$.
Then $Y={\mathbb R}^n\times L$, where $L$ is a group topologically isomorphic to
the character group of the group $K$. Denote elements of the group $Y$ by $y=(s, l)$,
where $s\in{\mathbb R}^n$, $l\in L$. In so doing,
\begin{equation}\label{04.09.6}
\hat\lambda_j(0, l)=
\begin{cases}
1, & \text{\ if\ }\   l=0,
\\  0, & \text{\ if\ }\ l\neq 0, \ \ j=1, 2.
\end{cases}
\end{equation}

Denote by $\alpha_K$ the restriction of   $\alpha$ to the subgroup $K$ and consider equation (\ref{2a}) for the characteristic functions
$\hat\lambda_j(s, l)$ on the group $Y={\mathbb R}^n\times L$. We have
\begin{equation}\label{25.03.nn1}
\hat\lambda_1(s_1+s_2, l_1+l_2)\hat\lambda_2(s_1+\tilde\alpha_{{\mathbb R}^n}s_2, l_1+\tilde\alpha_K l_2)$$$$=
\hat\lambda_1(s_1-s_2, l_1-l_2)\hat\lambda_2((s_1-\tilde\alpha_{{\mathbb R}^n}s_2, l_1-\tilde\alpha_K l_2), \ \  s_j\in {\mathbb R}^n, \ \ l_j \in L.
\end{equation}
Substitute $s_1=s_2=s$, $l_1=l$, $l_2=-l$ in
(\ref{25.03.nn1}). We obtain
\begin{equation}\label{04.09.7}
\hat\lambda_1(2s, 0)\hat\lambda_2((I+\tilde\alpha_{{\mathbb R}^n})s, (I-\tilde\alpha_K)l)=
\hat\lambda_1(0, 2l)\hat\lambda_2((I-\tilde\alpha_{{\mathbb R}^n})s, (I+\tilde\alpha_K)l), \ \  s \in {\mathbb R}^n, \ \ l \in L.
\end{equation}
 Assume that  $l\ne 0$. Then $2l\neq 0$, and (\ref{04.09.6}) implies that
 $\hat\lambda_1(0, 2l)=0$. Hence, it follows from  (\ref{04.09.7}) that
$\hat\lambda_1(2s, 0)\hat\lambda_2((I+\tilde\alpha_{{\mathbb R}^n})s, (I-\tilde\alpha_K)l)=0$.
Taking into account (\ref{04.09.2}), this implies that $$\hat\lambda_2((I+\tilde\alpha_{{\mathbb R}^n})s, (I-\tilde\alpha_K)l)=0.$$  Since $\alpha_{{\mathbb R}^n}$
 satisfies  condition (\ref{1}) and $\alpha({\mathbb R}^n)={\mathbb R}^n$, we have $I+\alpha_{{\mathbb R}^n}\in {\rm Aut}({{\mathbb R}^n})$. Hence, $I+\tilde\alpha_{{\mathbb R}^n}\in {\rm Aut}({{\mathbb R}^n})$ and
\begin{equation}\label{04.09.8x}
\hat\lambda_2(s, (I-\tilde\alpha_K)l)=0, \ \ s\in {{\mathbb R}^n}, \ \ l\in L, \ \ l\neq 0.
\end{equation}
In particular,
\begin{equation}\label{04.09.8}
\hat\lambda_2(0, (I-\tilde\alpha_K)l)=0, \ \  l\in L, \ \ l\neq 0.
\end{equation}
It follows from (\ref{04.09.6}) and (\ref{04.09.8}) that
\begin{equation}\label{04.09.9}
{\rm Ker}(I-\tilde\alpha_K)=\{0\}.
\end{equation}
Since $K$ is a compact subgroup of the group $G$,  the group  $K$ is topologically isomorphic to a compact totally disconnected Abelian group
   of the form $(\ref{n9})$.
  This implies that the group $L$ is   topologically isomorphic to a group of the form
 (\ref{m1}). Then, it follows from (\ref{04.09.9}) that $I-\tilde\alpha_K\in{\rm Aut}(L)$.
 Hence, we find from (\ref{04.09.8x}) that
\begin{equation}\label{04.09.10}
\hat\lambda_2(s, l)=0, \ \ s\in {{\mathbb R}^n}, \ \ l\in L, \ \ l\neq 0.
\end{equation}
Similarly, we obtain
\begin{equation}\label{04.09.11}
\hat\lambda_1(s, l)=0, \ s\in {{\mathbb R}^n}, \ \ l\in L, \ \ l\neq 0.
\end{equation}
 It follows from (\ref{04.09.2}), (\ref{04.09.10}), and (\ref{04.09.11}) that
 \begin{equation}\label{04.09.12}
\hat\lambda_j(s, l)=
\begin{cases}
\exp\{i\langle t_j, s\rangle-\langle A_js,
s\rangle\}, & \text{\ if\ }\   l=0,
\\  0, & \text{\ if\ }\ l\neq 0, \ \ j=1, 2.
\end{cases}
\end{equation}
Taking into account (\ref{04.09.12}), and returning from the distributions $\lambda_j$
to the distributions $\mu_j$,
we prove the statement of the lemma. $\blacksquare$

\textbf{\textit{Proof of Theorem \ref{th2}}}     We have $Y={\mathbb R}^n\times L\times H$, where $L$ and $H$ are groups
topologically isomorphic to the character groups of the groups $D$ and $G$ respectively. Let
$$G=\mathop{\mbox{\rm\bf P}}\limits_{p\in {\cal P}, \ p>2}{G}_{p},
$$
where $G_p$  is a finite $p$-group. Then
$$H=\mathop{\mbox{\rm\bf P}^*}\limits_{p\in {\cal P}, \ p>2}{H}_{p},
$$
where $H_p$ is a group
topologically isomorphic to the character group  of the group  $G_p$.
 Put $$A_m=\mathop{\mbox{\rm\bf P}}\limits_{p\in {\cal P}, \ 2<p\le m}{G}_{p}, \ \ B_m=\mathop{\mbox{\rm\bf P}}\limits_{p\in {\cal P}, \ 2<p\le m}{H}_{p}, \ \
 C_m=\mathop{\mbox{\rm\bf P}}\limits_{p\in {\cal P}, \ p>m}{G}_{p},
 \ \ m=3, 4, \dots$$
Since ${\mathbb R}^n$ is the connected component of zero of the group $X$, we have $\alpha({\mathbb R}^n)={\mathbb R}^n$. It follows from the definition of the product topology that $\alpha(C_m)\subset G$ for some $m$. Taking into account that $C_m$ is the only subgroup of $G$ which is topologically isomorphic to $C_m$, we have  $\alpha(C_m)=C_m$, and hence
$\alpha({\mathbb R}^n\times C_m)={\mathbb R}^n\times C_m$.

Represent the group $X$ in the form $X={\mathbb R}^n\times D\times A_m\times C_m$. Consider the factor-group
$X/({\mathbb R}^n\times C_m)$. Denote by $[x]$ its elements.
  Note that   the factor-group $X/({\mathbb R}^n\times C_m)$
  is topologically isomorphic to the group  $D\times A_m$. Since
  $\alpha({\mathbb R}^n\times C_m)={\mathbb R}^n\times C_m$, the topological automorphism
 $\alpha$ induces a topological automorphism $[\alpha]$ on the factor-group
 $X/({\mathbb R}^n\times C_m)$ by the formula
  $[\alpha] [x]=[\alpha x]$. It is easy to see that any continuous
monomorphism of the group ${\mathbb R}^n\times C_m$ to itself is a topological
automorphism. This easily implies that
\begin{equation}\label{2d}
{\rm Ker}(I+[\alpha])=\{0\}.
\end{equation}
 The character group of the factor-group $X/({\mathbb R}^n\times C_m)$ is
 topologically isomorphic to the annihilator
$A(Y, {\mathbb R}^n\times C_m)=L\times B_m$. Since $\alpha({\mathbb R}^n\times C_m)={\mathbb R}^n\times C_m$, we have
$\tilde\alpha(L\times B_m)=(L\times B_m)$.
Note also that the restriction of
 $\tilde\alpha$ to the subgroup $L\times B_m$ is the adjoint automorphism to $[\alpha]$.

Let $\eta_1$ and $\eta_2$ be independent random variables with values
in the factor-group  $X/({\mathbb R}^n\times C_m)$ such that
the characteristic functions of their distributions coincide with
the restrictions to the subgroup $L\times B_m$ of the characteristic functions $\hat\mu_j(y)$.   By Lemma \ref{le1}, the conditional distribution
 of the linear form $M_2 = \eta_1 + [\alpha]\eta_2$ given
  $M_1 = \eta_1 + \eta_2$ is symmetric. Obviously, the factor-group
  $X/({\mathbb R}^n\times C_m)$  is discrete and contains no elements of
  order 2.
Hence, taking into account (\ref{2d}) and applying  Theorem B to the discrete group
 $X/({\mathbb R}^n\times C_m)$, we obtain that
the restrictions of the characteristic functions $\hat\mu_j(y)$
to the subgroup $L\times B_m$ are of the form
$$
\hat\mu_j(y)=
\begin{cases}
(g_j, y), & \text{\ if\ }\   y\in E,
\\  0, & \text{\ if\ }\ y\notin E, \ \ j=1, 2,
\end{cases}
$$
where $g_j\in D\times A_m$, $E$ is an open subgroup in $L\times B_m$, which
is the annihilator of
a finite subgroup in the group $X/({\mathbb R}^n\times C_m)$. Moreover,
\begin{equation}\label{06.09.3}
\tilde\alpha(E)=E.
\end{equation}
Using the same reasoning as in the proof of  Lemma \ref{le8},    we can replace the independent random variables
$\xi_j$ by $\eta_j=\xi_j-x_j$  with distributions $\lambda_j$ in such a way that   the conditional distribution of the linear form
$N_2 = \eta_1 + \alpha\eta_2$
given $N_1 = \eta_1 + \eta_2$  is symmetric,   and
the restrictions of the characteristic functions
$\hat\lambda_j(y)$ to the subgroup $L\times B_m$ are of the form
\begin{equation}\label{06.09.2}
\hat\lambda_j(y)=
\begin{cases}
1, & \text{\ if\ }\   y\in E,
\\  0, & \text{\ if\ }\ y\notin E, \ \ j=1, 2.
\end{cases}
\end{equation}
By Lemma \ref{le7}, (\ref{06.09.2}) implies that $\sigma(\lambda_j)\subset A(X, E)$.
It is easy to see that $A(X, E)={\mathbb R}^n\times F$, where $F$ is
a compact totally disconnected Abelian group
   of the form $(\ref{n9})$. It follows from (\ref{06.09.3}) that $\alpha({\mathbb R}^n\times F)={\mathbb R}^n\times F$. Applying  Lemma \ref{le8} to the group ${\mathbb R}^n\times F$, we get the required  representation for
   the distributions $\lambda_j$, and hence, for the distributions $\mu_j$ too. $\blacksquare$

 An analogue of   the Heyde  theorem is valid not only for locally compact
 Abelian groups of the form (\ref{09.09.15.1}).
The following statement holds true (compare with  \cite[\S 3]{FG}).
\begin{proposition} \label{pr2} There exists a locally compact Abelian group
$X$ containing no elements of order $2$   such that $X$   is not topologically isomorphic to any group of the form
$(\ref{09.09.15.1})$,  and   the following statement holds. If $\alpha$ is
       a topological automorphism of
the group  $X$ satisfying   condition $(\ref{1})$, and
$\xi_1$ and  $\xi_2$ are independent random variables with values in
       $X$  and distributions $\mu_1$ and $\mu_2$,  then
the symmetry of the conditional distribution of the linear form
$L_2 = \xi_1 + \alpha\xi_2$ given  $L_1 = \xi_1 +
\xi_2$ implies that
$\mu_j=m_K*E_{x_j}$, where $K$ is a compact subgroup
of   $X$,
$x_j\in X$, $j=1, 2$. Moreover, $\alpha (K)=K$.
\end{proposition}
\textbf{\textit{Proof}}    Denote by
 $p_n$
the $n$th prime number. Put $T=\mathop{\mbox{\rm\bf
P}}\limits_{n=2}^\infty{\mathbb Z}(p_n^2)$  and consider $T$ only as "algebraic" group. We do  not consider   $T$ as a topological group equipped with the product topology.

Consider a subgroup   $K_n=\{0, p_n, 2p_n, \dots, (p_n-1)p_n\}$ of ${\mathbb Z}(p^2_n)$.  It is obvious that $K_n$ is isomorphic to
${\mathbb Z}(p_n)$.
Denote by
$Y$ a subgroup of $T$ consisting of  all elements  $y=(y_2, \dots, y_n, \dots)$, $y_n\in{\mathbb Z}(p^2_n)$, such that
$y_n\notin K_n$  only for a finite number indexes
$n$. Put
$$
F=\mathop{\mbox{\rm\bf P}}\limits_{n=2}^\infty K_n.
$$
Then $F$ is a subgroup of $Y$. Consider the group $F$ in the product topology.
Obviously, $F$ is a compact Abelian group. Define a topology on $Y$ as follows. A subset $U$ in $Y$ is called open if  for each $y\in U$ there is  a neighborhood  of zero $V_y$ in $F$ such that $y+V_y\subset U$.
Then $Y$ is a second countable locally compact Abelian group, and $F$
is an open subgroup in $Y$.
It is obvious that $Y$ is a non-discrete, non-compact and totally disconnected group. Denote by
$X$   the character group of the group $Y$.

Since multiplication by 2 is a topological isomorphism of the group $Y$, the group $X$ contains no elements of order 2.  To see that the group $X$ is not topologically isomorphic to a group of the form   (\ref{09.09.15.1}),  suppose  the contrary.  It is obvious that then in (\ref{09.09.15.1}) $n=0$, and then
   the group $Y$ can be represented in the form
$Y=L\times H$,
where the group $L$ is compact, and the group $H$ is discrete.
  It is easy to see that there exists $m$, such that $\mathop{\mbox{\rm\bf P}}\limits_{n=m}^\infty K_n\subset L$. This implies that the group $H$ contains only a finite number different of subgroups  $K_n$, and hence, $H$ is a finite group.
  It follows from this that $Y$ is a compact group. The obtained contradiction shows that the group $X$ is not topologically isomorphic to a group of the form  $(\ref{09.09.15.1})$.

We shall verify that
\begin{equation}\label{10.09.15.2}
{(I+\tilde\alpha)(F)}=F.
\end{equation}
Note that $\alpha$  satisfies  condition $(\ref{1})$ if and only if
\begin{equation}\label{10.09.15.1}
\overline{(I+\tilde\alpha)(Y)}=Y.
\end{equation}
Put $B=\mathop{\mbox{\rm\bf P}^*}\limits_{n=2}^\infty{\mathbb
Z}(p^2_n)$, and consider $B$ as a subgroup of the group $Y$ in the topology induced on $B$ by   the topology of $Y$.  Denote by $\tilde\alpha_B$ the restriction of $\tilde\alpha$ to the subgroup $B$.
Assume that ${\rm Ker}(I+\tilde\alpha_B)\ne\{0\}$. This implies that
$(I+\tilde\alpha_B)(B)$ is a proper subgroup of $B$. It follows from this that $(I+\tilde\alpha)(B)$ is contained in a proper closed subgroup of the group $Y$. Since $B$ is a dense subgroup in $Y$,
this implies that $(I+\tilde\alpha)(Y)$ is also contained in a proper closed subgroup of the group $Y$, contrary to
$(\ref{10.09.15.1})$. Hence, ${\rm Ker}(I+\tilde\alpha_B)=\{0\}$,
and then $I+\tilde\alpha_B$ is an automorphism of the group $B$. Put $C=\mathop{\mbox{\rm\bf P}^*}\limits_{n=2}^\infty K_n$, and consider $C$ as a subgroup of the group $F$ in the topology induced on $C$ by   the topology of $Y$.   Then $C$ is a dense subgroup in $F$. Since $I+\tilde\alpha_B$ is an automorphism of the group $B$, we have
\begin{equation}\label{10.09.15.5}
(I+\tilde\alpha)(C)=C.
\end{equation}
Take $y_0\in F$. The subgroup $C$ is dense in $F$. Hence, there exists a sequence  $ y_j\in C$ such that $  y_j\rightarrow  y_0$. It follows from (\ref{10.09.15.5}) that
$y_j=(I+\tilde\alpha)z_j$, where $z_j\in C$. Since $F$ is a compact group, there exist a sequence $\{j_k\}$ and an element $z\in F$ such that
$z_{j_k}\rightarrow z$. Hence $ y_0=(I+\tilde\alpha)  z$.
Thus, (\ref{10.09.15.2}) is proved.

By Lemma \ref{le1}, the characteristic functions
  $\hat\mu_j(y)$ satisfy equation
(\ref{2a}). Put
$\nu_j = \mu_j
* \bar \mu_j$. Then
 $\hat \nu_j(y) = |\hat \mu_j(y)|^2 \ge 0$  for all $y \in Y$, and
 the characteristic functions  $\hat \nu_j(y)$
also satisfy equation (\ref{2a}). It is easy to see that Proposition \ref{pr2} will be proved if we prove
that $\nu_1=\nu_2=m_K$, where $K$ is a compact subgroup of the group
 $X$, and $\alpha(K)=K$.
Therefore we can solve equation
(\ref{2a}), assuming that  $\hat \mu_j(y) \ge 0$ for all  $y \in Y$. Denote by $\tilde\alpha_F$ the restriction of $\tilde\alpha$ to the subgroup $F$. It is obvious that $\tilde\alpha_F\in {\rm Aut}(F)$.  Consider the restriction of equation (\ref{2a}) to the subgroup $F$. We have
\begin{equation}\label{25.03.2}
\hat\mu_1(u+v )\hat\mu_2(u+\tilde\alpha_F v )=
\hat\mu_1(u-v )\hat\mu_2(u-\tilde\alpha_F v), \ \ u, v \in F.
\end{equation}
The group $F$ is topologically isomorphic to the character group of the discrete group
$\mathop{\mbox{\rm\bf P}^*}\limits_{n=2}^\infty K_n,$ and   $\tilde\alpha_F$   is the adjoint
automorphism   to an automorphism
$\bar\alpha$ of the group $\mathop{\mbox{\rm\bf P}^*}\limits_{n=2}^\infty K_n.$ Note that $(\ref{10.09.15.2})$ implies that  $\bar\alpha$ satisfies  condition
 (\ref{1}). Taking into account Lemma \ref{le1}, by    applying Theorem B to
 the discrete group
$\mathop{\mbox{\rm\bf P}^*}\limits_{n=2}^\infty K_n,$ we conclude from (\ref{25.03.2}) that the restrictions of the characteristic functions
$\hat\mu_j(y)$ to the subgroup $F$ are of the form
$$
\hat\mu_j(y)=
\begin{cases}
1, & \text{\ if\ }\   y\in E,
\\  0, & \text{\ if\ }\ y\notin E, \ \ j=1, 2,
\end{cases}
$$
where $E$ is an open subgroup of $F$. In so doing,
\begin{equation}\label{10.09.15.4}
\tilde\alpha_F(E)=E.
\end{equation}
Put $M=A(X, E)$. It follows from (\ref{10.09.15.4}) that $\alpha(M)=M$.
By Lemma \ref{le7}, $\sigma(\mu_j)\subset M$, $j=1, 2$. Thus we reduced the proof of Proposition \ref{pr2} to the case when the random variables
$\xi_j$
take values in the group $M$.  The character group of the group $M$
is topologically isomorphic to the factor-group $Y/E$. Since $E$ is an open subgroup of $F$, there  is some natural number $m$ such that the group $Y/E$ is  topologically isomorphic to a group of the form
$$
\mathop{\mbox{\rm\bf P}}\limits_{n=2}^m{\mathbb
Z}(p^{l_n}_n)\times\mathop{\mbox{\rm\bf P}^*}\limits_{n=m+1}^\infty{\mathbb
Z}(p_n),
$$
where $l_n\in\{0, 1, 2\}$. Hence, the group $M$ is topologically isomorphic to a group of the form (\ref{n9}). The proposition   follows now from Theorem \ref{th1}. $\blacksquare$

\section{  The Heyde  theorem fails on compact connected
Abelian groups}

   We will prove in this section that there does not exist a compact connected
Abelian group on which even a weak analogue of   the Heyde   theorem holds (see below Theorem \ref{th3}).

Let $n$ be a natural number. Denote by $f_n:X \rightarrow X$ an endomorphism of
the group $X$, defined by the formula  $f_nx=nx$ for all $x\in X$.
Put $f_n(X)=X^{(n)}$. Denote by $\mathbb{Q}$ the additive group of rational numbers
considering in the discrete topology. Let $\text{\boldmath $a$}=(a_0,a_1,\dots,a_n,\dots)$ be an
arbitrary infinite sequence of integers, where each of $a_n$ is
greater than  1.
Denote by  $\Sigma_\text{\boldmath $a$}$  the corresponding ${\text{\boldmath $a$}}$-adic solenoid (\!\!\cite[(10.12)]{Hewitt-Ross}). The group $\Sigma_\text{\boldmath $a$}$  is compact,  connected  and has dimension 1 (\!\!\cite[(10.13), (24.28)]{Hewitt-Ross}).
Denote by $H_{\text{\boldmath $a$}}$
a subgroup of $\mathbb{Q}$ of the form
$$\displaystyle{H_{\text{\boldmath $a$}}=\left\{{m\over
a_0a_1\cdots a_n}: \ n=0, 1, \dots; \ m=0, \pm 1, \pm 2, \dots\right\}}.
$$
The character group of the group $\Sigma_{\text{\boldmath $a$}}$ is topologically isomorphic to the group $H_{\text{\boldmath $a$}}$ (\!\!\cite[(10.12)]{Hewitt-Ross}). Note that if $\text{\boldmath $a$}=(2,3,4,\dots)$, then $H_{\text{\boldmath $a$}}=\mathbb{Q}$.

The main result of this section is   the following statement.
\begin{theorem}\label{th3}  Let $X$ be a compact connected Abelian group containing no elements of order $2$.  Then there exist a topological automorphism $\alpha$ of the group $X$ satisfying   condition $(\ref{1})$ and
   independent random variables $\xi_1$ and  $\xi_2$  with values in
       $X$  and distributions $\mu_1$ and $\mu_2$ such that
 the conditional distribution of the linear form
$L_2 = \xi_1 + \alpha\xi_2$ given  $L_1 = \xi_1 +
\xi_2$  is symmetric, whereas  $\mu_j\notin
\Gamma(X)*I(X)$, $j=1, 2$.
\end{theorem}
We note that if $\mu\in
\Gamma(X)*I(X)$, then $\mu$ is invariant with respect to
a compact subgroup $K$ of the group $X$ and, under the natural homomorphism
  $X\rightarrow X/K$, induces a Gaussian distribution on the factor-group  $X/K$. Thus, Theorem \ref{th3} shows that even a weak analogue the  Heyde   theorem is not valid for compact connected Abelian groups.

\textbf{\textit{Proof of Theorem \ref{th3}}}   Since $X$ is a compact connected
Abelian group,  $Y$ is a discrete torsion free Abelian group.
Two cases are possible:  either  the group $X$ contains   elements of finite order
or $X$ is a torsion-free group.

{1}. Assume that the group $X$ contains    elements of finite order. Denote by $p$ the minimum of orders of nonzero elements of $X$.  Then $p$ is a prime and   taking into account that $X$ contains no elements of order 2, $p\ge 3$. Since  $X$ is a connected Abelian group, we have
$X^{(n)} = X$
 for all natural $n$. This implies that if  ${\rm Ker}f_n=\{0\}$, then
 $f_n \in {\rm Aut}(X)$. Set $\alpha=-f_{p-1}$. Then  $\alpha\in {\rm Aut}(X)$. Since $I+\alpha=-f_{p-2}$, we have $I+\alpha\in {\rm Aut}(X)$, so that $\alpha$ satisfies  condition $(\ref{1})$.
It follows from $Y^{(p)}=A(Y, {\rm Ker}f_p)$ that $Y^{(p)}$ is a proper
subgroup in the group  $Y$.
Take $y_0\in Y$ such that
$y_0\notin Y^{(p)}$. Since $p \ge 3$, the numbers 2 and $p$ are mutually prime.  Therefore there exist integers $m$  and $n$ such that $2m+pn=1$.
It follows from this that $y_0=2my_0+pny_0$.
 This implies that   $2y_0\notin Y^{(p)}$.

Consider on the group $X$ the function
$$
\rho(x)=1+{\rm Re}(x, y_0), \ \ x\in X.
$$
Then $\rho(x)\ge 0$ for all $x\in X,$ and
$$
\int\limits_X\rho(x)dm_X(x)=1.
$$
Denote by $\mu$ the distribution on $X$ with the density
$\rho(x)$ with respect to $m_X$.
It is easy to see that the characteristic function
$\hat\mu(y)$ is of the form
\begin{equation}\label{08a1}
\hat \mu(y)=
\begin{cases}
1, & \text{\ if\ }\   y=0,
\\  {1\over 2}, & \text{\ if\ }\ y=\pm y_0,
\\ 0,& \text{\ if\ }\ y\not\in\{0, \pm y_0\}.
\end{cases}
\end{equation}
It is obvious that $\mu\not\in
\Gamma(X)*I(X)$.

Let $\xi_1$ and  $\xi_2$ be independent identically distributed
   random variables   with values in the group
       $X$  and distribution $\mu$.  Let us verify that
 the conditional distribution of the linear form
$L_2 = \xi_1 + \alpha\xi_2$ given  $L_1 = \xi_1 +
\xi_2$  is symmetric.  By Lemma \ref{le1}, it suffices to verify that the characteristic function $\hat\mu(y)$
satisfies  equation (\ref{08a2}).
Obviously, equation (\ref{08a2}) holds true if  $v=0$. So, assume that $v\ne 0$. Then the left-hand side of equation (\ref{08a2}) vanishes.
Otherwise, taking into account (\ref{08a1}), we get
$u+v\in \{0, \pm y_0\}$ and $u+\alpha v\in \{0, \pm y_0\}$. This implies that $(I-\tilde\alpha) v=f_pv\in \{0, \pm y_0, \pm 2y_0\}$, but that's impossible because it contradicts $y_0, 2y_0\notin Y^{(p)}$. Reasoning similarly we obtain that if $v\ne 0$, then the right-hand side of equation (\ref{08a2}) vanishes too. Thus, in the case when the group $X$ contains    elements of finite order,  the theorem
is proved.

 {2}. Assume that  $X$ is a torsion-free group. Taking into account the structure theorem for compact Abelian torsion-free group (\!\!\cite[(25.8)]{Hewitt-Ross}) and the fact that   $X$ is  connected, we get that for some $\mathfrak{n}$ the group $X$     is topologically isomorphic to a group of the form
$(\Sigma_{\text{\boldmath $a$}})^\mathfrak{n}$, where ${\text{\boldmath $a$}}=(2, 3, 4, \dots).$
Obviously, arguing as in the proof of necessity  in  Theorem \ref{th1}  it suffices to prove the theorem for the group $
X=\Sigma_{\text{\boldmath $a$}}$, where ${\text{\boldmath $a$}}=(2, 3, 4, \dots)$. Then the group $Y$ is topologically isomorphic to the group $\mathbb{Q}$. In order not to complicate the notation we assume that $Y=\mathbb{Q}$.

Let $H$ be a subgroup of the group $Y$ of the form
$$\displaystyle{H=\left\{{m
\over 5^n}: \ n=0, 1, \dots; \ m=0, \pm 1, \pm 2, \dots\right\}}.
$$
 Consider on the group $H$   the function
\begin{equation}
\label{08a3} g(h) =
\begin{cases}
1, & \text{\ if\ }\ \ h \in H^{(2)},
\\ c, & \text{\ if\ }\ \ h \notin H^{(2)},
\end{cases}
\end{equation}
where $-1 < c < 1,$ $c\ne 0$.  It is easy to verify that $g(h)$ is a positive definite function. Consider on the group $Y$ the function
\begin{equation}
\label{08a4} f(y) =
\begin{cases}
  g(y), & \text{\ if\ }\ \ y \in H,
\\ 0, & \text{\ if\ }\ \ y \notin  H.
\end{cases}
\end{equation}
Then $f(y)$ is also a positive definite function
 (\!\!\cite[(32.43)]{Hewitt-Ross2}). By the Bochner theorem, there exists
 a distribution $\mu \in{\rm
M}^1(X)$ such that $\hat\mu(y) = f(y).$ It is obvious that $\mu
\notin \Gamma(X)*I(X).$

Put $\alpha=-f_{2}^{-2}$. Then $\alpha\in {\rm Aut}(X)$, and we  have
$I+\alpha=f_3f_{2}^{-2}$.
Obviously, $I+\alpha\in {\rm Aut}(X)$, so that
$\alpha$ satisfies  condition  $(\ref{1})$.

Let $\xi_1$ and  $\xi_2$ be independent identically distributed
   random variables   with values in the group
       $X$  and distribution $\mu$.  Verify that
 the conditional distribution of the linear form
$L_2 = \xi_1 + \alpha\xi_2$ given  $L_1 = \xi_1 +
\xi_2$  is symmetric. This is equivalent to the fact that
the conditional distribution of the linear form
$N_2 = 4\xi_1 -\xi_2$
given $N_1 = \xi_1 + \xi_2$ is symmetric. By Lemma \ref{le1}, it suffices to verify that the characteristic  function $\hat\mu(y)$
satisfies the equation (\ref{08a2}) which takes the form
\begin{equation}
\label{08a5} f(u + 4v)f(u - v) =
f(u - 4v)f(u + v),
\ \ u, v \in Y.
\end{equation}
It follows from  (\ref{08a3}) and (\ref{08a4}) that for $u, v \in H$ the equation (\ref{08a5}) becomes an equality.
Note that if either $u \in H$, $v \notin H$ or $v \in H$, $u \notin H,$ then
$u\pm v\notin H$, and  it follows from (\ref{08a4}) that both sides of equation (\ref{08a5}) are equal to zero.  So, equation (\ref{08a5}) becomes also an equality.
Assume that $u, v \notin H.$ Then the left-hand side in (\ref{08a5}) is equal to zero. In the opposite case (\ref{08a4}) implies that
 $u + 4v,  u - v \in
H.$ It follows from this that $5v \in H$, and hence, $v \in H$.
The obtained contradiction shows that the left-hand side in (\ref{08a5}) is equal to zero. Arguing similarly, we verify that
if $u, v \notin H,$ then the right-hand side in (\ref{08a5}) is also equal to zero. $\blacksquare$

In connection with Theorem \ref{th3}, it is interesting to note that the following statement is true (see \cite[Theorem 3]{Fe7}).

{\it  Let  $X$ be a locally compact Abelian group containing
no elements of order $2$.  Let $\alpha$ be a topological automorphism of the group
$X$ satisfying   condition $(\ref{1})$.
Let  $\xi_1$ and  $\xi_2$ be independent random variables with values in
$X$   and distributions   $\mu_1$ and $\mu_2$ with non-vanishing characteristic functions.
  If  the conditional distribution of the linear form
$L_2 = \xi_1 + \alpha\xi_2$ given  $L_1 = \xi_1 +
\xi_2$  is symmetric, then $\mu_j\in
\Gamma(X)$, $j=1, 2$.}

We will complement Theorem \ref{th3} with two propositions. The first proposition shows that,
generally speaking, we can not replace in Theorem \ref{th3}
 the statement:

 {\it {\rm (A)} \ Then there exist a topological automorphism $\alpha$ of the group $X$ satisfying   condition $(\ref{1})$ and
   independent random variables $\xi_1$ and  $\xi_2$  with values in
       $X$  and distributions $\mu_1$ and $\mu_2$ such that
 the conditional distribution of the linear form
$L_2 = \xi_1 + \alpha\xi_2$ given  $L_1 = \xi_1 +
\xi_2$  is symmetric, whereas  $\mu_j\notin
\Gamma(X)*I(X)$, $j=1, 2$.}

by the statement:

 {\it   {\rm (B)}  \ Then for each topological automorphism  $\alpha$ of the group $X$, $\alpha\ne I$, satisfying   condition $(\ref{1})$ there exist
   independent random variables $\xi_1$ and  $\xi_2$  with values in
       $X$  and distributions $\mu_1$ and $\mu_2$ such that
 the conditional distribution of the linear form
$L_2 = \xi_1 + \alpha\xi_2$ given  $L_1 = \xi_1 +
\xi_2$  is symmetric, whereas  $\mu_j\notin
\Gamma(X)*I(X)$, $j=1, 2$.}

The second proposition shows that there is a compact connected Abelian group $X$ containing no elements of order $2$ when we can  replace in Theorem \ref{th3} (A) by (B). We should exclude the case when $\alpha=I$, because if a  compact Abelian group $X$ contains no elements of order 2, then ${Y^{(2)}}=Y$, and by Lemma \ref{le1}, the symmetry of the conditional distribution
of the linear form  $L_2 = \xi_1 + \xi_2$ given
 $L_1 = \xi_1 + \xi_2$ implies that  $\xi_1$ and $\xi_2$ have degenerate
 distributions.
\begin{proposition} \label{pr3} There exist a compact connected Abelian group $X$ containing no elements of order $2$ and a topological automorphism $\alpha$  of the group $X$, $\alpha\not=I$, satisfying  condition $(\ref{1})$ such that the following statement holds: If $\xi_1$ and  $\xi_2$ are
   independent random variables    with values in
       $X$  and distributions $\mu_1$ and $\mu_2$ such that
 the conditional distribution of the linear form
$L_2 = \xi_1 + \alpha\xi_2$ given  $L_1 = \xi_1 +
\xi_2$  is symmetric, then $\mu_j\in
\Gamma(X)*I(X)$, $j=1, 2$.
\end{proposition}
To prove Proposition \ref{pr3} we need two lemmas.
\begin{lemma} [\!\!{\protect\cite[Lemma 6]{My1}}] \label{le9}
Let  $X$ be a locally compact Abelian group, and let $\alpha$
be a topological automorphism of the group $X$.
Let  $\xi_1$ and  $\xi_2$ be independent random variables with values in
  $X$.
   If the conditional distribution of the linear form    $L_2 = \xi_1 +
\alpha\xi_2$  given $L_1 = \xi_1 +
\xi_2$  is symmetric,   then the linear forms
$M_1=(I+\alpha)\xi_1+2\alpha\xi_2$ and
$M_2=2\xi_1+(I+\alpha)\xi_2$  are independent.
\end{lemma}
\begin{lemma}    [\!\!{\protect\cite {Fe8}, see also \cite[Theorem 7.10]{Fe5a}}]
\label{le10}
Let  $X$ be a locally compact Abelian group. Assume that the connected component of zero of the
group $X$
 contains no elements of order $2$. Let   $\xi_1$ and $\xi_2$ be
 independent random variables
with values in  $X$ and  distributions $\mu_1$ and $\mu_2$. If the sum $\xi_1+\xi_2$ and the difference $\xi_1-\xi_2$ are independent, then
$\mu_j\in \Gamma(X)*I(X)$, $j=1, 2$.
\end{lemma}
\textbf{\textit{Proof of Proposition \ref{pr3}}}   Let $X=\Sigma_{\text{\boldmath $a$}}^2$, where
${\text{\boldmath $a$}}=(2, 3, 4, \dots)$. Since the character group of
the group $X$ is topologically isomorphic to the group  $\mathbb{Q}^2$, we have
${f_n\in {\rm Aut}(X)}$
for any natural $n$. This implies in particular, that the group
 $X$ contains no elements of order 2. Denote by $x=(a, b)$, where $a, b\in \Sigma_{\text{\boldmath $a$}}$, elements of the group $X$.
Define the mapping $\alpha: X\rightarrow X$  by the formula
$$\alpha(a, b)=(-3a+4b, 2a-3b), \ (a, b)\in X.$$
Then $\alpha$ is a topological automorphism of the group $X$.   Since $I+\alpha\in {\rm Aut}(X)$,   condition (\ref{1}) is fulfilled. By Lemma \ref{le9},
the symmetry of the conditional distribution of the linear form    $L_2 = \xi_1 +
\alpha\xi_2$  given $L_1 = \xi_1 +
\xi_2$ implies that the linear forms
$M_1=(I+\alpha)\xi_1+2\alpha\xi_2$ and
$M_2=2\xi_1+(I+\alpha)\xi_2$  are independent.
 Note that coefficients of the linear forms  $M_1$ and  $M_2$ are topological automorphisms of the group $X$.
Consider the new independent random variables $\eta_1=(I+\alpha)\xi_1$ and
$\eta_2=2\alpha\xi_2$. The independence of the linear forms
$M_1$ and $M_2$ implies that  the linear forms $N_1=\eta_1+\eta_2$
and $N_2=2(I+\alpha)^{-1}\eta_1+f_2^{-1}\alpha^{-1}(I+\alpha)\eta_2$ are independent.
Hence,  the linear forms $P_1=\eta_1+\eta_2$ and $P_2=\eta_1+f_2^{-2}\alpha^{-1}(I+\alpha)^2\eta_2$ are also independent. It is easy to verify that
 $f_{2}^{-2}\alpha^{-1}(I+\alpha)^2=-I$. Hence, $P_2=\eta_1-\eta_2$.
Denote by
  $\lambda_j$ the distribution of the random variable $\eta_j$.
  Since $X$ is a connected Abelian
group
 containing no elements of order $2$, by Lemma \ref{le10},  $\lambda_j\in
\Gamma(X)*I(X)$, $j=1, 2$. In view of $I+\alpha, 2\alpha\in {\rm Aut}(X)$, it follows from $\lambda_1=(I+\alpha)\mu_1$ and $\lambda_2=2\alpha\mu_2$     that $\mu_j\in
\Gamma(X)*I(X)$, $j=1, 2$. $\blacksquare$
\begin{proposition} \label{pr4} There exists a compact connected Abelian group $X$  containing no elements of order $2$ such that for each topological automorphism   $\alpha$
    of  the group  $X$, $\alpha\ne I$,  satisfying  condition $(\ref{1})$, there exist independent random variables $\xi_1$ and  $\xi_2$  with values in
       $X$  and distributions $\mu_1$ and  $\mu_2$ such that
 the conditional distribution of the linear form
$L_2 = \xi_1 + \alpha\xi_2$ given  $L_1 = \xi_1 +
\xi_2$  is symmetric, while $\mu_j\notin
\Gamma(X)*I(X)$, $j=1, 2$.
\end{proposition}

\textbf{\textit{Proof}}   Let $X=\Sigma_\text{\boldmath $a$}$, where
 $\text{\boldmath $a$}=(2,2,2,\dots)$.
The group  $Y$ is topologically isomorphic to a subgroup of the group  $\mathbb{Q}$ of the form $$\displaystyle H_{\text{\boldmath $a$}}=\left\{{m
\over 2^n}: \ n=0, 1, \dots; \ m=0, \pm 1, \pm 2, \dots\right\}.$$
Since $Y^{(2)}=Y$, the group $X$  contains no elements of order 2.
 It is easy to see that ${\rm Aut}(X)={\{\pm f_{  2^k}, \pm f_{  2^k}^{-1}: k=0, 1, 2, \dots\}.}$ Obviously, there exist only two topological automorphisms
 $\alpha$, $\alpha\ne I$, of the group $X$ satisfying   condition (\ref{1}). They are
 $\alpha=-f_{2}$ and $\alpha=-f_{2}^{-1}$. Let $\alpha=-f_{2}$. It is easy to see that the group $X$ contains an element of order 2. Then, as has been proved in the case    1  of Theorem  \ref{th3},
 there exist independent random variables $\xi_1$ and  $\xi_2$  with values in $X$  and distributions $\mu_1$ and  $\mu_2$ such that
 the conditional distribution of the linear form
$L_2 = \xi_1 -2\xi_2$ given  $L_1 = \xi_1 +
\xi_2$  is symmetric, while $\mu_j\notin
\Gamma(X)*I(X)$, $j=1, 2$.

The proof of the existence of independent random variables $\xi_1$ and  $\xi_2$  with values in $X$  and distributions $\mu_1$ and  $\mu_2$ such that
 the conditional distribution of the linear form
$L_2 = \xi_1 -f_{2}^{-1}\xi_2$ given  $L_1 = \xi_1 +
\xi_2$   is symmetric, while $\mu_j\notin
\Gamma(X)*I(X)$, $j=1, 2$,  is reduced to the previous case. $\blacksquare$

\begin{remark}\label{r2}
The Heyde  theorem is true for some compact connected Abelian groups if  we additionally assume that the characteristic functions of random variables do not vanish. In particular, the following statement was proved   in  \cite{Fe2020}.

\begin{theorema}  Let $X=\Sigma_\text{\boldmath $a$}$ be   an \text{\boldmath $a$}-adic solenoid   containing no elements
of order $2$,  and let  $\alpha$ be a topological automorphism  of the group  $X$ satisfying   condition $(\ref{1})$.
Let $\xi_1$ and $\xi_2$ be independent random variables with values in the group
       $X$  and distributions
  $\mu_1$ and $\mu_2$ with nonvanishing characteristic functions. If the conditional distribution of the linear form
$L_2 = \xi_1 + \alpha\xi_2$ given  $L_1 = \xi_1 +
\xi_2$  is symmetric, then
  $\mu_j\in \Gamma(X)$,  $j=1, 2$.
\end{theorema}
\end{remark}

\noindent{\bf Acknowledgments}        My special thanks  go to the  referee for   careful  reading of my article and many useful suggestions and  remarks. I appreciate it very much.

\vskip 1 cm

\noindent G. M. Feldman

\medskip

\noindent e-mail: feldman@ilt.kharkov.ua

\medskip

\noindent B. Verkin Institute for Low Temperature \\
Physics and Engineering of the\\
National Academy of Sciences of Ukraine\\
47, Nauky Ave  \\
Kharkiv, 61103  \\
Ukraine

\end{document}